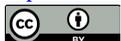
Scientific
Research
Publishing

# A New Extrapolation Economy Cascadic Multigrid Method for Image Restoration Problems

## Zhaoteng Chu, Ziqi Yan, Chenliang Li

School of Mathematics and Computational Science, Guilin University of Electronics Technology, Guilin, China
Email: chenli@guet.edu.cn

## Abstract

In this paper, a new extrapolation economy cascadic multigrid method is proposed to solve the image restoration model. The new method combines the new extrapolation formula and quadratic interpolation to design a nonlinear prolongation operator, which provides more accurate initial values for the fine grid level. An edge preserving denoising operator is constructed to remove noise and preserve image edges. The local smoothing operator reduces the influence of staircase effect. The experiment results show that the new method not only improves the computational efficiency but also ensures good recovery quality.

## 1. Introduction

In the process of image formation, transmission, recording and processing, it is easy to be affected by equipment, environment, human factors and so on. It is very important to recover the noisy and blurred image. Therefore, it is very meaningful to construct a fast and effective image restoration algorithm, which has become one of the research hotspots in the field of image processing.

Consider the following two-dimensional image degradation model,

$$g(x, y) = h(x, y) * u(x, y) + \eta(x, y), \tag{1}$$

where $h(x, y)$ is an impulse response function, which is called the point spread function in the spatial domain, $*$ represents the two-dimensional convolution operation, $u(x, y)$ represents the input image to be processed, $\eta(x, y)$ represents the additive noise, and $g(x, y)$ represents the degraded image.

Ignoring the noise $\eta(x, y)$, Equation (1) can be expressed as

$$g(x, y) = h(x, y) * u(x, y).\tag{2}$$

In practical problems, it is usually necessary to restore the blur- and noise-contaminated image. Therefore, we usually consider the following problem,

$$h(x, y) * u(x, y) = g^{\delta}(x, y),\tag{3}$$

where $g^{\delta}(x, y)$ represents the blur- and noise-contaminated image.

Equation (3) is an ill-posed problem. In the case of direct solution, a small disturbance of the right-hand side may lead to great changes, so the recovered image is usually of low quality. The Tikhonov regularization method can eliminate this instability by solving the following model

$$\min_{u} \frac{1}{2} \int_{-\infty}^{+\infty} \int_{-\infty}^{+\infty} \left( h * u - g^{\delta} \right)^2 \mathrm{d}x\mathrm{d}y + \alpha \psi(u),\tag{4}$$

where the first term is the fidelity term, $\psi(u)$ is the regularization term, also known as constraints, and $\alpha$ is the regularization parameter, which is used to balance the fidelity term and regularization term. The discrete regularization term is denoted by $D(u)$, and the selection of the regularization operator $D(u)$ is crucial for the quality of image restoration. The most representative one is the nonlinear anisotropic diffusion equation proposed by Perona and Malik in [1], and the regularization operator is expressed as follows

$$D(u) = div \left[ c \left( |\nabla u| \right) \nabla u \right],\tag{5}$$

the diffusion coefficients are chosen as follows,

$$c \left( |\nabla u| \right) = \frac{1}{1 + \left( \dfrac{|\nabla u|}{k} \right)^2},\tag{6}$$

where $k$ is the threshold.

For nonlinear Equation (4), [2] and [3] propose the explicit and semi-implicit discretization schemes, which have good preserving effect on the edges of the image and can provide high quality restoration of images. However, whether it is the explicit or semi-implicit discretization schemes its computation cost is very large. In addition, the selection of the regularization parameters and the time intervals are not easy.

Equation (2) is discretized to obtain the following linear system

$$g = Hu, \ H \in \mathbb{R}^{n \times n}, \ u, g \in \mathbb{R}^{n},\tag{7}$$

where $H$ represents the fuzzy operator, which is usually a real symmetric block Teoplitiz matrix, $u$ is the original image, and $g$ represents the blur and noise-free image.

The blur- and noise-contaminated image $g$ is obtained by adding random normally distributed noise $e$ with mean 0 and standard deviation 1, such as

$$g^{\delta} = g + e. \tag{8}$$

Thus, Equation (3) is discretized to obtain the following linear system

$$Hu = g^{\delta}, \ H \in \mathbb{R}^{n \times n}, \ u, g^{\delta} \in \mathbb{R}^{n}. \tag{9}$$

we set the solution is $u^{\delta}$.

For the linear Equation (9), the matrix $H$ is usually a large-scale singular matrix. The system of linear equations is usually an ill-posed problem. Although we can solve Equation (9) by some iterative method, due to the loss of high-frequency information, these iterative methods may produce ringing effect and cannot accurately recover edges, which seriously affects the quality of image restoration.

In order to solve the problems lots of effective methods have been proposed.

In [4], a standard iterative method is proposed to compute several approximate solutions and then extract a new candidate solution from the linear subspace expanded by these approximate solutions. This method can be used for large-scale problems. In [5], LSQR algorithm and RRGMRES algorithm are proposed to solve the large-scale linear ill-posed problems, and corresponding stopping criteria are proposed, which can obtain meaningful approximate solutions. In [6], the regularization properties of the global GMRES method are analyzed, and a regularized global GMRES method is proposed for solving linear discrete ill-posed problems with multiple right-hand sides contaminated by errors. In [7], NL-means algorithm is proposed to correct the ladder effect.

Multigrid method is one of the most effective methods for solving elliptic boundary value problems. The cascadic multigrid method is proposed in [8]. The advantage of this method lies in its simplicity and efficiency. In [9] presents an economical cascadic multigrid method. Compared with the conventional cascadic multigrid method, the proposed method can get the same efficient solutions with fewer iterations. [10] [11] [12] have proposed some extrapolation cascadic multigrid methods. The methods use the new extrapolation with quadratic interpolation, which can provide a better initial value for the fine grid level.

Based on the advantages of these methods, many scholars have applied them to image restoration. [13] [14] have proposed multigrid method for image restoration. The step effect is reduced by introducing wavelet soft threshold denoising and smoothing. In [15], a cascadic multilevel method is proposed to solve a linear discrete ill-posed problem with error contamination on the right side. Numerical experiments show that the algorithm is effective in image restoration. In [16], a class of nonlinear edge-preserving continuation operators are proposed to be used as prolongation operators for cascadic multigrid methods. This method can reduce the ringing effect and restore good image edges. In [17], the selection of fuzzy operator and the smooth method are further studied, and it is shown that the cascadic multilevel method has good efficiency for image restoration.

Combining the advantages of the new extrapolation method and the cascadic multigrid method, a new extrapolation economy cascadic multigrid method is proposed for image restoration in this paper. The experimental results show that this method not only improves computational efficiency but also ensures good recovery quality. The new method has the following advantages:

1) The new extrapolation formula and the quadratic interpolation are combined to provide more accurate values for the fine mesh level.

2) Fewer iterations are costed on each grid level without affecting the accuracy.

3) The edge preserving denoising operator can preserve the edge of image while denoising.

4) The local smoothing operator reduces the ladder effect and ensures good recovery quality.

The structure of this paper is as follows. In Section 2, introduces the basic framework of the new extrapolation economy cascadic multigrid. In Section 3, the linear prolongation operator, the new extrapolation formula prolongation operator combined with quadratic interpolation and the smoothing operator generated by the weighted least squares approximation is introduced. In Section 4, the edge preserving denoising operator is introduced. In Section 5, some numerical results are given to illustrate the effectiveness of the new method. In Section 6, a summary is presented.

## 2. New Extrapolation Economy Cascadic Multigrid Method

### Selecting a Template

Let $V = \left[ v^{(1)}, v^{(2)}, \cdots, v^{(m)} \right]^{\mathrm{T}}$, we define the following norm

$$\|V\| = \left( \frac{1}{m} \sum_{i=1}^{m} \left( v^{(i)} \right)^2 \right)^{1/2}, \tag{10}$$

For the Equation (9), we can solve it by using the conjugate gradient method (CG) or the minimal residual method (MR). Let the initial iterative value $u_0^l = 0$, here $l$ represents the grid level. We use the following stopping rule to terminate the iteration.

Let $\delta = \|e\|$, for a given $\delta \geq 0$, and $c > 1$ is a constant and independent of $\delta$, when $\left\| g^{\delta} - Hu_k^l \right\| \leq c\delta$, stop the iteration.

By using CG or MR to solve the linear model directly, although it cannot effectively restore the edge which seriously affects the quality of image restoration, it is simple and the time is short. In view of these two advantages, we can use CG or MR as the smoothing operator in the new extrapolation economic cascadic multigrid method.

For the sake of discussion, we set $u$ as an $n^2 \times 1$ column vector and $N = n^2$, $G_1 \subset \cdots \subset G_l \subset \cdots \subset G_L$ is a nested subspace. The subscript $l$ indicates the number of the grid level, $G_1$ denotes the coarsest grid level, and $G_L$ indicates the finest grid level. Accordingly, we have $N_1 < \cdots < N_l < \cdots < N_L = N$. When $N$ is odd, $N_{l-1} = \left( \sqrt{N_l} + 1 \right)^2 / 4$, when $N$ is even, $N_{l-1} = N_l / 4$.

In the cascadic multigrid algorithm, the following linear models need to be solved iteratively at each grid level

$$H_l u_l = g_l^\delta, \tag{11}$$

where $H_l$ is the representation of the fuzzy operator $H$ in the $G_l$.

Set $\left(g_l^\delta\right)^{(k)}$ be the $k^{\text{th}}$ element of $g_l^\delta$. Let

$$\left(g_l^\delta\right)^{(k)} = \left(g_{l+1}^\delta\right)^{(2k-1)}, 1 \le k \le N_l \tag{12}$$

We define the restriction operator $R_l$, which satisfies

$$g_l^\delta = R_l g_{l+1}^\delta. \tag{13}$$

When $l = L$, $g_l^\delta = g^\delta$.

$H_l$ is generated in the following way,

$$H_l = R_l H_{l+1} R_l^*. \tag{14}$$

where $R_l^*$ is the adjoint of $R_l$, when $l = L$, $H_l = H$.

In [16], a cascadic multiresolution method is proposed. Numerical examples show that the method has good performance in image restoration. It is explained that for Equation (9), lots of iterations will affect the recovery effect. Therefore, the number of iterations is particularly critical, and the cascadic multigrid method proposed by [9] is of reference significance. We combine the advantages of the algorithms of [9] and [16] propose the following improved economic cascadic multigrid method.

**Algorithm 1**: Improved Economical Cascadic Multigrid Method for Image Restoration (IECMG)

Step 1: Set $u_1^0 = 0$, $u_1^\delta = C_1^{m_1} u_1^0$.

Step 2: When $l = 2, 3, \cdots, L$, do:

1) Prolongation: $\tilde{u}_l^0 = L_l u_{l-1}^\delta$ or $\tilde{u}_l^0 = P_l u_{l-1}^\delta$.

2) Local smoothing: $u_l^S = S_l \tilde{u}_l^0$.

3) Preserving edge and removing noise: $u_l^{DS} = D_l S_l \tilde{u}_l^0$.

4) Smoothing: $u_l^0 = u_l^{DS}$, $u_l^\delta = C_l^{m_l} u_l^0$.

Step 3: $u^\delta = u_L^\delta$.

In [9] [10] [11] [12], the economic cascadic multigrid method and the new extrapolation cascadic multigrid method are proposed respectively. The new extrapolation interpolation method can provide better initial values for the next level and has better convergence. Therefore, we propose the following new extrapolation economy cascadic multigrid method.

**Algorithm 2**: New Extrapolation Economy Cascadic Multigrid method for Image Restoration (EECMG)

Step 1: Set $u_l^0 = 0$, $u_l^\delta = C_l^{m_l} u_l^0$, $l = 1, 2$.

Step 2: When $l = 2, 3, \cdots, L$, do:

1) New extrapolation: $\tilde{u}_l^0 = \Pi_l \left(u_{l-1}^\delta, u_l^\delta\right)$.

2) Quadratic interpolation prolongation: $\tilde{u}_{l+1}^0 = P_l \tilde{u}_l^0$.

3) Local smoothing: $u_{l+1}^S = S_{l+1} \tilde{u}_{l+1}^0$.

4) Preserving edge and removing noise: $u_{l+1}^{DS} = D_{l+1} S_{l+1} \tilde{u}_{l+1}^0$.

5) Smoothing: Let $u_{l+1}^0 = u_{l+1}^{DS}$, $u_{l+1} = C_{l+1}^{m_{l+1}} u_{l+1}^0$

Step 3: $u^\delta = u_L^\delta$

The number of iterations $m_l$ on each level is selected as follows [9].

1) If $l > L_0$, then $m_l = \left\lceil m_0 \left( L - L_0 \right)^2 \beta^{L-l} \right\rceil$.

2) If $l \le L_0$, then $m_l = \left\lceil m_*^{1/2} \left( L - (2 - \varepsilon_0) l \right) h_l^{-2} \right\rceil$

Where $l = 1, 2, \cdots, L$, $L_0 = \dfrac{1}{2} L$, $m_0 = 1$, $\beta = 4$, $m_* = 1$, $\varepsilon_0 = \dfrac{1}{2}$,

$h_l = \left( \dfrac{1}{2} \right)^{l-1}$, $[t]$ represents the smallest positive integer not less than $t$.

Note: In the above algorithms, $C_l^{m_l}$ denotes smoothing $m_l$ times by using CG or MR. In Sections 3 and Sections 4, we will introduce the linear prolongation operator $L_l$, nonlinear prolongation operator $P_l$, new extrapolation operator $\Pi_l$, local smoothing operator $S_l$ and edge preserving denoising operator $D_l$ in Algorithm 1 and Algorithm 2 in detail.

# 3. Prolongation Operator and Local Smoothing Operator

## 3.1. Linear Prolongation Operator

Let $u_{(i,j)}^l$ denote the pixel value in $(i, j)$ on the level *l*. The linear operator $L_l$ can be defined as follows,

$$
\begin{aligned}
u_{(i,j)}^l &= u_{((i+1)/2,(j+1)/2)}^{l-1}, && \text{when } i, j \text{ both odd;} \\
u_{(i,j)}^l &= \frac{1}{2} \left( u_{((i+1)/2,j/2)}^{l-1} + u_{((i+1)/2,j/2+1)}^{l-1} \right), && \text{when } i \text{ is odd, } j \text{ is even;} \\
u_{(i,j)}^l &= \frac{1}{2} \left( u_{(i/2,(j+1)/2)}^{l-1} + u_{(i/2+1,(j+1)/2)}^{l-1} \right), && \text{when } i \text{ is even, } j \text{ is odd;} \\
u_{(i,j)}^l &= \frac{1}{2} \left( u_{(i/2,j/2)}^{l-1} + u_{(i/2+1,j/2+1)}^{l-1} \right), && \text{when } i, j \text{ are both even.}
\end{aligned}
\tag{15}
$$

Because the accuracy of the image on the fine grid level by piecewise linear interpolation is not high, a nonlinear prolongation operator, combined a new extrapolation formula with quadratic interpolation, is constructed to solve this problem.

## 3.2. New Extrapolation Quadratic Interpolation Prolongation Operator

Let us take the one-dimensional triplet grid $G_l (l = 1, 2, 3)$ as an example, the pixel is denoted by $p_i^l$, and the corresponding pixel value is denoted by $u_i^l$. Let $G_1$ be the coarsest grid level and contains pixels $I_1 = \left( p_i^1, p_{i+1}^1 \right)$. The pixels contained in $G_2$ and $G_3$ are represented as follows,

$$
I_2 = \left( p_i^2, p_{(i+1)/2}^2, p_{i+1}^2 \right); \quad I_3 = \left( p_i^3, p_{(i+1)/4}^3, p_{(i+1)/2}^3, p_{(i+3)/4}^3, p_{i+1}^3 \right),
\tag{16}
$$

the corresponding pixel value are expressed as,

$$
\begin{aligned}
&G_1 = \left( u_i^1, u_{i+1}^1 \right); \quad G_2 = \left( u_i^2, u_{(i+1)/2}^2, u_{i+1}^2 \right); \\
&G_3 = \left( u_i^3, u_{(i+1)/4}^3, u_{(i+1)/2}^3, u_{(i+3)/4}^3, u_{i+1}^3 \right).
\end{aligned}
\tag{17}
$$

We combine the new extrapolation formula with quadratic interpolation, which is called the new extrapolation quadratic interpolation method. This can provide better initial values on the fine grid level $G_3$. The details are as follows (see **Figure 1**).

1) New extrapolation $u_3 = \Pi_2(u_1, u_2)$,

$$u_i^3 = \frac{1}{4}\left(5u_i^2 - u_i^1\right);$$

$$u_{i+1}^3 = \frac{1}{4}\left(5u_{i+1}^2 - u_{i+1}^1\right);$$

$$u_{i+\frac{1}{2}}^3 = u_{i+\frac{1}{2}}^2 + \frac{1}{8}\left[\left(u_i^2 - u_i^1\right) + \left(u_{i+1}^2 - u_{i+1}^1\right)\right].$$

2) Quadratic interpolation $u_3 = P_3(u_1, u_2, u_3)$,

$$u_{i+\frac{1}{4}}^3 = \frac{3}{8}u_i^3 + \frac{3}{4}u_{i+\frac{1}{2}}^3 - \frac{1}{8}u_{i+1}^3$$

$$= \frac{1}{16}\left[\left(9u_i^2 + 12u_{i+\frac{1}{2}}^2 - u_{i+1}^2\right) - \left(3u_i^1 + u_{i+1}^1\right)\right];$$

$$u_{i+3/4}^3 = -\frac{1}{8}u_i^3 + \frac{3}{4}u_{i+1/2}^3 + \frac{3}{8}u_{i+1}^3$$

$$= \frac{1}{16}\left[\left(9u_{i+1}^2 + 12u_{i+1/2}^2 - u_i^2\right) - \left(3u_{i+1}^1 + u_i^1\right)\right].$$

Next, we apply this idea to the two-dimension problems. As shown in **Figure 2**, the hollow circles represent the pixels on the $G_{I-2}$, the solid black dots represent the pixel values on the $G_{I-1}$ grid level, and the rectangular boxes represent the pixel values on the $G_I$ grid level. We use the new extrapolation formula to get the image value at the black node and the quadratic interpolation to obtian the image value at the rectangular box.

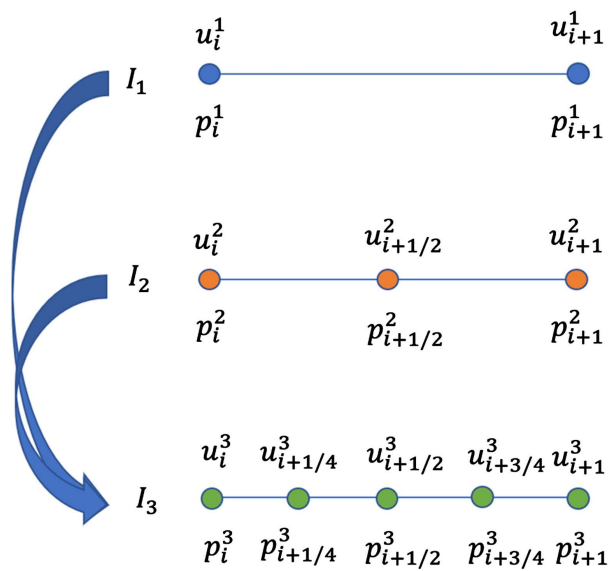

**Figure 1.** New extrapolation quadratic interpolation method for the 1-dimensional.

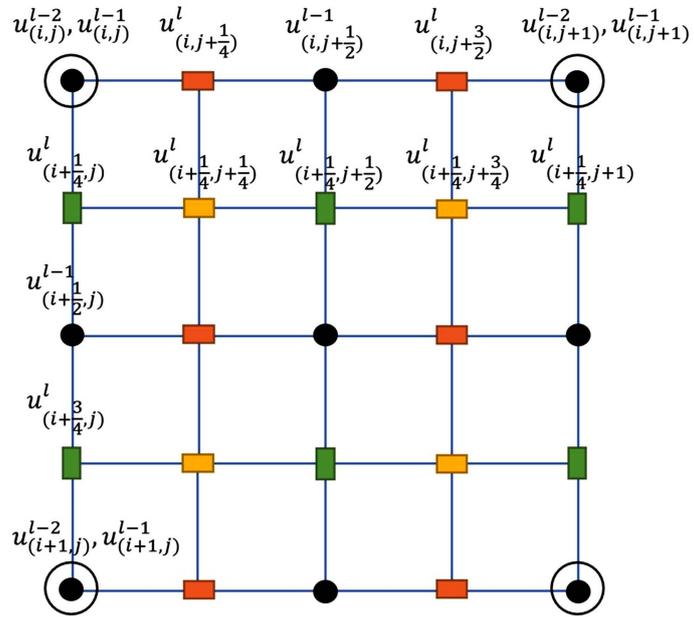

**Figure 2.** New extrapolation quadratic interpolation method for the 2-dimensional.

1) New extrapolation $u_l = \Pi_l\left(u_{l-1}, u_{l-2}\right)$,

$$u^l_{(i,j)} = \frac{1}{4}\left(5u^{l-1}_{(i,j)} - u^{l-2}_{(i,j)}\right);$$

$$u^l_{(i,j+1)} = \frac{1}{4}\left(5u^{l-1}_{(i,j+1)} - u^{l-2}_{(i,j+1)}\right);$$

$$u^l_{\left(i,j+\frac{1}{2}\right)} = u^{l-1}_{\left(i,j+\frac{1}{2}\right)} + \frac{1}{8}\left[\left(u^{l-1}_{(i,j)} - u^{l-2}_{(i,j)}\right) + \left(u^{l-1}_{(i,j+1)} - u^{l-2}_{(i,j+1)}\right)\right];$$

$$u^l_{(i+1/2,j)} = u^{l-1}_{(i+1/2,j)} + \frac{1}{8}\left[\left(u^{l-1}_{(i,j)} - u^{l-2}_{(i,j)}\right) + \left(u^{l-1}_{(i+1,j)} - u^{l-2}_{(i+1,j)}\right)\right].$$

2) Quadratic interpolation $u_l = P_l\left(u_l, u_{l-1}, u_{l-2}\right)$,

$$u^l_{\left(i,j+\frac{1}{4}\right)} = \frac{1}{16}\left[\left(9u^{l-1}_{(i,j)} + 12u^l_{\left(i,j+\frac{1}{2}\right)} - u^{l-1}_{(i,j+1)}\right) - \left(3u^{l-2}_{(i,j)} + u^{l-2}_{(i,j+1)}\right)\right];$$

$$u^l_{\left(i,j+\frac{3}{4}\right)} = \frac{1}{16}\left[\left(9u^{l-1}_{(i,j+1)} + 12u^l_{\left(i,j+\frac{1}{2}\right)} - u^{l-1}_{(i,j)}\right) - \left(3u^{l-2}_{(i,j+1)} + u^{l-2}_{(i,j)}\right)\right];$$

$$u^l_{\left(i+\frac{1}{4},j\right)} = \frac{1}{16}\left[\left(9u^{l-1}_{(i,j)} + 12u^l_{\left(i+\frac{1}{2},j\right)} - u^{l-1}_{(i+1,j)}\right) - \left(3u^{l-2}_{(i,j)} + u^{l-2}_{(i+1,j)}\right)\right];$$

$$u^l_{\left(i+\frac{3}{4},j\right)} = \frac{1}{16}\left[\left(9u^{l-1}_{(i+1,j)} + 12u^l_{\left(i+\frac{1}{2},j\right)} - u^{l-1}_{(i,j)}\right) - \left(3u^{l-2}_{(i+1,j)} + u^{l-2}_{(i,j)}\right)\right];$$

$$u^l_{(i+1/4,j+1/4)} = \frac{3}{8}u^l_{(i+1/4,j)} + \frac{3}{4}u^l_{(i+1/4,j+1/2)} - \frac{1}{8}u^l_{(i+1/4,j+1)};$$

$$u^l_{(i+1/4,j+3/4)} = -\frac{1}{8}u^l_{(i+1/4,j)} + \frac{3}{4}u^l_{(i+1/4,j+1/2)} + \frac{3}{8}u^l_{(i+1/4,j+1)}.$$

### 3.3. Local Smoothing Operator

In order to get a smoother image, the local weighted least squares approximation is used to smooth the image in [2]. As shown in **Figure 3**, the pixel value $u_{(i,j)}^l$ is related to the image value of the surrounding eight points. The process is represented by the operator $S_l$ as follows.

First, we introduce the weighting function,

$$\omega_{(i,j)}^l(s,t) = \exp\left(-\left(u_{(i,j)}^l - u_{(i+s,j+t)}^l\right)^2\right), \ \ s,t \in \{0,\pm 1\}, \tag{18}$$

then we solve the locally weighted least squares problem,

$$\min_{a_0,a_1,a_2} \sum_{s,t\in\{0,\pm 1\}} \left(u_{(i+s,j+t)}^l - \left(a_0 + a_1 s + a_2 t\right)\right)^2 \omega_{(i,j)}^l(s,t). \tag{19}$$

Set $\hat{a}_0, \hat{a}_1, \hat{a}_2$ represent the solutions of the above problem, then the updated pixel value is $u_{(i,j)}^l = \hat{a}_0$.

## 4. Edge Preserving Denoising Operator

The isotropic diffusion equation diffuses equally along the tangential direction and the normal direction of the image edge, so it cannot protect the edge, nor can it preserve the original fine structure of the image well, so that the image becomes blurred. In order to deal with these shortcomings, an anisotropic diffusion equation is proposed. One is the nonlinear anisotropic diffusion equation proposed by Perona and Malik, abbreviated as P-M equation [1], which is as follows,

$$\begin{cases} \dfrac{\partial u}{\partial t} = \mathrm{div}\left[g\left(|\nabla u|\right)\nabla u\right], \\ u(x,y,t)\big|_{t=0} = u_0(x,y), \end{cases} \tag{20}$$

where the $u(x,y,t)$ is the image of the $t$ moment, $g\left(|\nabla u|\right) \in [0,1]$ is the diffusion coefficient. The more classical diffusion coefficient is as follows,

$$g\left(|\nabla u|\right) = \frac{1}{1 + \left(\dfrac{|\nabla u|}{k}\right)^2} \tag{21}$$

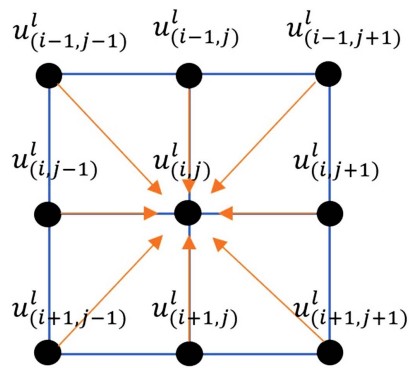

**Figure 3.** Local smoothing method.

Set $c(x,y) = g(|\nabla u|)$, discretization of (20) is as follows,

$$\frac{\partial u}{\partial t} = \text{div}\left[g(|\nabla u|)\nabla u\right] = \nabla\left[c(x,y)\nabla u\right] = \left(c(x,y)u_y\right)_x + \left(c(x,y)u_y\right)_y$$

$$= \frac{\partial}{\partial x}\left[c(x,y)\frac{u\left(x+\frac{\Delta x}{2},y\right) - u\left(x-\frac{\Delta x}{2},y\right)}{\Delta x}\right]$$

$$+ \frac{\partial}{\partial y}\left[c(x,y)\frac{u\left(x,\frac{\Delta x}{2}+y\right) - u\left(x,y+\frac{\Delta x}{2}\right)}{\Delta y}\right]$$

$$= \frac{1}{\Delta x}\left[c\left(x+\frac{\Delta x}{2},y\right)\frac{1}{\Delta x}\left(u(x+\Delta x,y) - u(x,y)\right)\right.$$

$$\left. - c\left(x-\frac{\Delta x}{2},y\right)\frac{1}{\Delta x}\left(u(x,y) - u(x-\Delta x,y)\right)\right]$$

$$+ \frac{1}{\Delta y}\left[c\left(x,\frac{\Delta x}{2}+y\right)\frac{1}{\Delta y}\left(u(x,y+\Delta y) - u(x,y)\right)\right.$$

$$\left. - c\left(x,y-\frac{\Delta x}{2}\right)\frac{1}{\Delta y}\left(u(x,y) - u(x,y-\Delta y)\right)\right], \qquad (22)$$

Taking $\Delta x = 1, \Delta y = 1$, we have,

$$\frac{u_{(i,j)}^{n+1} - u_{(i,j)}^{n}}{\tau} = \frac{c_{(i+1,j)}^{n} + c_{(i,j)}^{n}}{2}\left(u_{(i+1,j)}^{n} - u_{(i,j)}^{n}\right) + \frac{c_{(i-1,j)}^{n} + c_{(i,j)}^{n}}{2}\left(u_{(i-1,j)}^{n} - u_{(i,j)}^{n}\right)$$

$$+ \frac{c_{(i,j+1)}^{n} + c_{(i,j)}^{n}}{2}\left(u_{(i,j+1)}^{n} - u_{(i,j)}^{n}\right) + \frac{c_{(i,j-1)}^{n} + c_{(i,j)}^{n}}{2}\left(u_{(i,j-1)}^{n} - u_{(i,j)}^{n}\right) \qquad (23)$$

$$= \sum_{(k,l)\in Q(i,j)} \frac{c_{(k,l)}^{n} + c_{(i,j)}^{n}}{2}\left(u_{(k,l)}^{n} - u_{(i,j)}^{n}\right),$$

where $Q(i,j)$ denotes the set of four neighboring points centered on $(i,j)$. So according to the above equation we have,

$$u_{(i,j)}^{n+1} = u_{(i,j)}^{n} + \tau\left[\frac{c_{(i+1,j)}^{n} + c_{(i,j)}^{n}}{2}u_{(i+1,j)}^{n} + \frac{c_{(i-1,j)}^{n} + c_{(i,j)}^{n}}{2}u_{(i-1,j)}^{n} + \frac{c_{(i,j+1)}^{n} + c_{(i,j)}^{n}}{2}u_{(i,j+1)}^{n}\right.$$

$$\left. + \frac{c_{(i,j-1)}^{n} + c_{(i,j)}^{n}}{2}u_{(i,j-1)}^{n} - \frac{c_{(i+1,j)}^{n} + c_{(i-1,j)}^{n} + c_{(i,j+1)}^{n} + c_{(i,j-1)}^{n} + 4c_{(i,j)}^{n}}{2}u_{(i,j)}^{n}\right]. \qquad (24)$$

Equation (24) can be expressed in matrix form as,

$$u^{n+1} = u^n + \tau A(u^n)u^n, \qquad (25)$$

where the elements of $A(u^n)$ are as follows,

$$a_{(i,j),(s,t)} = \begin{cases} c\dfrac{c_{(i,j)}^{n} + c_{(s,t)}^{n}}{2}, & (s,t) \in Q(i,j) \\[2mm] -\sum_{(m,n)\in Q(i,j)} \dfrac{c_{(i,j)}^{n} + c_{(m,n)}^{n}}{2}, & (s,t) = (i,j). \\[2mm] 0, & \text{other} \end{cases} \qquad (26)$$

# 5. Numerical Examples

We use numerical examples to illustrate the effectiveness of the new extrapolation economy cascadic multigrid method for restoring blur- and noise-contaminated image. The elements in the fuzzy operator $H$ are defined by the following function.

$$h_{(i,j)} = \begin{cases} \dfrac{1}{\sigma\sqrt{2\pi}}\exp\left(-\dfrac{(i-j)^2}{2\sigma^2}\right), & \text{if } |j-k| \leq \text{band} \\ 0, & \text{other} \end{cases}$$

$$h = \left[ h_{(i,j)} \right]_{i,j=1,\cdots,n},$$

$$H = h \otimes h.$$

where $\otimes$ denotes the Kronecker product. $\sigma$ represents the variance of the Gaussian point diffusion function. The "band" represents the half bandwidth of $h$. The degree of blurring is denoted by $b_i$, and the following cases are considered in this paper,

$$b_1 : \begin{cases} \sigma = 1 \\ \text{band} = 7 \end{cases}, \quad b_2 : \begin{cases} \sigma = 2 \\ \text{band} = 9 \end{cases}, \quad b_3 : \begin{cases} \sigma = 3 \\ \text{band} = 11 \end{cases}, \quad b_4 : \begin{cases} \sigma = 4 \\ \text{band} = 13 \end{cases}.$$

Defining the noise levels as follows,

$$\frac{\left\| u - u^\delta \right\|}{\|u\|} = \frac{\|e\|}{\|u\|} \approx v,$$

the noise level $v$ is set to the following values,

$$v_1 = 5 \times 10^{-2}, \quad v_2 = 1 \times 10^{-1}, \quad v_3 = 5 \times 10^{-1}, \quad v_4 = 8 \times 10^{-1}.$$

In this experiment, the blur and noise level of the image is denoted as $b_i v_i$ ($i = 1,2,3,4$). For example, $b_1 v_1$ denotes that $\sigma = 1$, band $= 7$, $v_1 = 5 \times 10^{-2}$.

The image restoration effect is measured by the following peak signal ratio,

$$\text{PSNR}(u,\hat{u}) = 20\log_{10}\frac{255}{\|u - \hat{u}\|}\text{dB}.$$

where $u$ represents blur- and noise-free image, $\hat{u}$ represents the recovered image.

Figure 4(a) is the original image, Figures 4(b)-(e) are the images with four different levels of blur and noise.

Firstly, CG algorithm and MR algorithm are used for image restoration. When the stopping rule is met, the iteration times of both algorithms are 4. The restoration effect of four different levels of blur- and noise-contaminated images are shown in Figure 5.

Table 1 and Figure 5 show that the MR algorithm has a faster computational speed than the CG algorithm in the same case. For $b_1 v_1$ and $b_2 v_2$, the recovery effect of CG and MR is almost the same. For $b_3 v_3$ and $b_4 v_4$, the MR algorithm has better recovery effect.

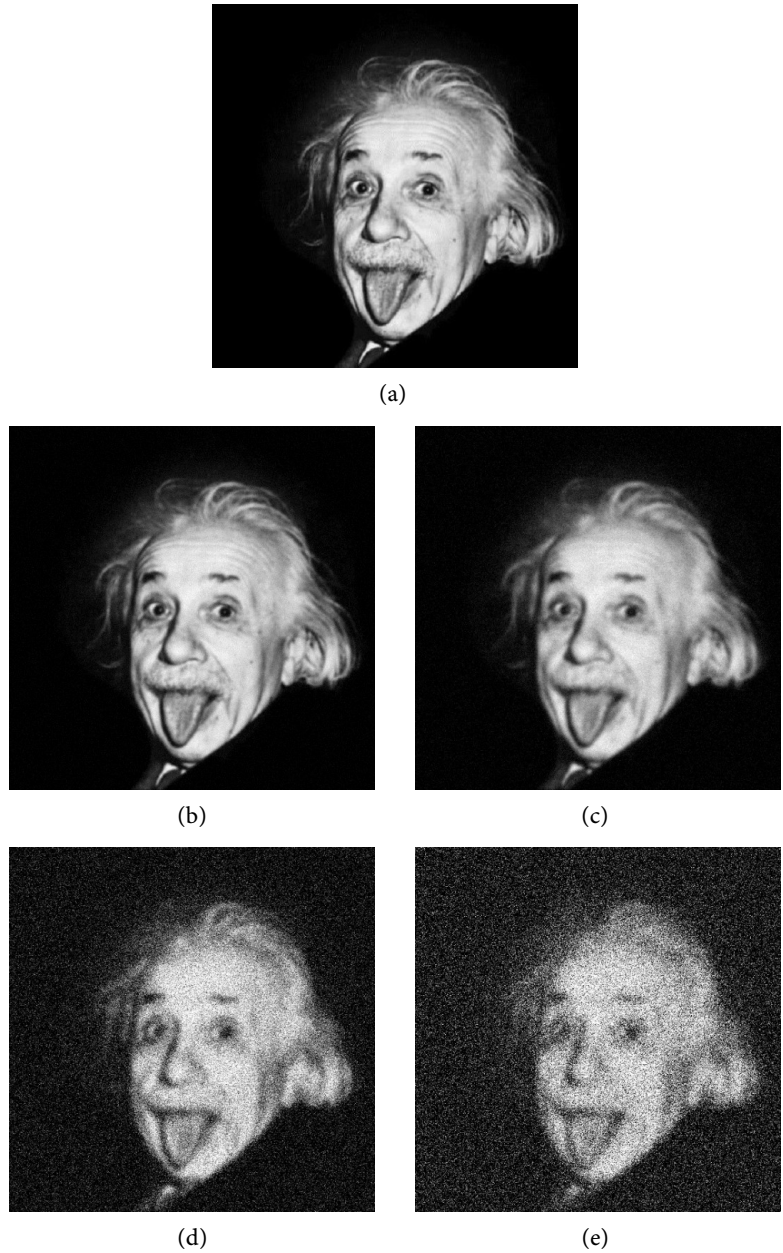

**Figure 4.** Original images and contaminated images with different levels of blur and noise. (a) Original image; (b) Blur and noise image $b_1v_1$; (c) Blur and noise image $b_2v_2$; (d) Blur and noise image $b_3v_3$; (e) Blur and noise image $b_4v_4$.

**Table 1.** Numerical results of CG and MR.

| Blur and noise levels | CG MR | | | |
| --- | --- | --- | --- | --- |
| | PSNR | Time | PSNR | Time |
| $b_1v_1$ | 31.5824 | 1.847 | 31.1796 | 0.996 |
| $b_2v_2$ | 26.1651 | 3.453 | 26.0547 | 1.406 |
| $b_3v_3$ | 17.1940 | 4.912 | 19.1158 | 1.994 |
| $b_4v_4$ | 15.1794 | 6.585 | 17.0716 | 2.553 |

In the following experiments, with different smoothers, Algorithm 1 and Algorithm 2 are abbreviated as shown in the Table 2. Algorithm 1 and Algorithm 2 with different smoother. Next, we compare the numerical results of restoration Algorithm 1 and Algorithm 2.

From Table 3 and Table 4, Figure 6 and Figure 7, it can be found that,

1) The quadratic interpolation can obtain better accuracy initial value. So IECMG-P is better than IECMG-L.

2) New extrapolation with quadratic interpolation can improve the accuracy of the initial values. So EECMG is better than IECMG.

3) With MR as smoother, EECMG (MR) and IECMG (MR) have the best restoration quality and spend the least recovery time.

Moreover, Figure 5, Figure 8 and Figure 9 show that EECMG (MR) is best, whether PSNR or CPU time.

Finally, we show the edge preserving ability of the new extrapolation cascadic multigrid method.

Figure 10(a) is the original image, Figure 10(b) is the third-level polluted image. Figure 10(c) is the restored image by using EECMG (MR), which PSNR = 27.9617. The good effect of edge preservation is shown in Figure 10(d).

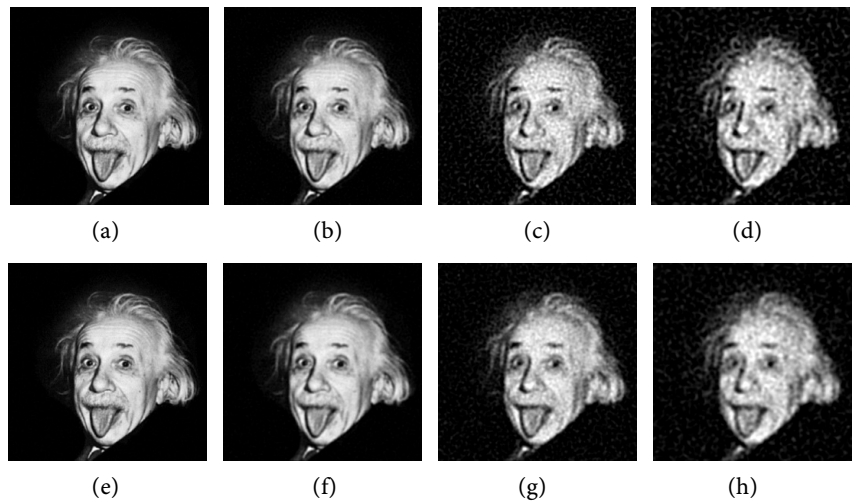

(a)          (b)          (c)          (d)

(e)          (f)          (g)          (h)

**Figure 5.** Numerical results of CG and MR. (a) CG, $b_1 v_1$; (*b*) CG, $b_2 v_2$; (c) CG, $b_3 v_3$; (d) CG, $b_4 v_4$; (e) MR, $b_1 v_1$; (f) MR, $b_2 v_2$; (g) MR, $b_2 v_2$; (h) MR, $b_4 v_4$.

**Table 2.** Algorithm 1 and Algorithm 2 with different smoother.

| | |
|---|---|
| IECMG-L (CG) | IECMG with liner interpolation and CG as smoother |
| IECMG-P (CG) | IECMG with quadratic interpolation and CG as smoother |
| EECMG (CG) | EECMG with new extrapolation quadratic interpolation and CG as smoother |
| IECMG-L (MR) | IECMG with liner interpolation and MR as smoother |
| IECMG-P (MR) | IECMG with quadratic interpolation and MR as smoother |
| EECMG (MR) | EECMG with new extrapolation quadratic interpolation and MR as smoother |

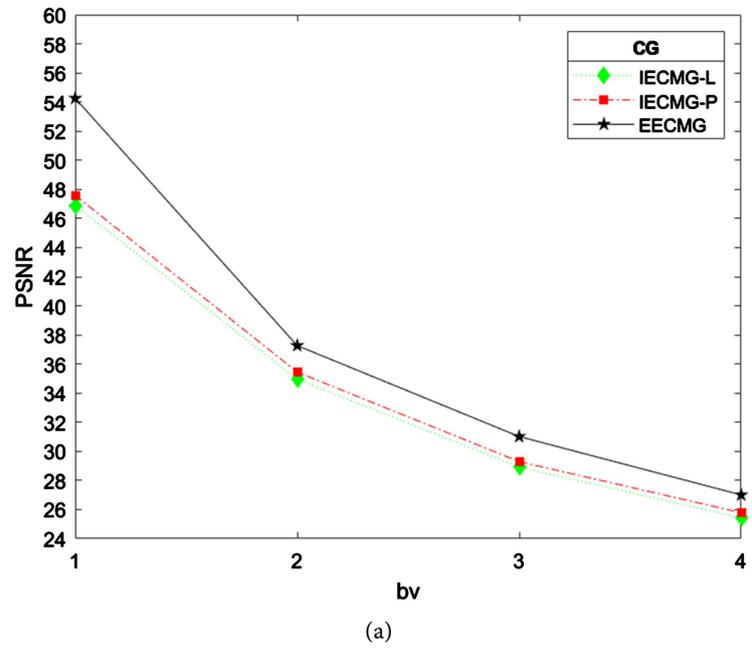

(a)

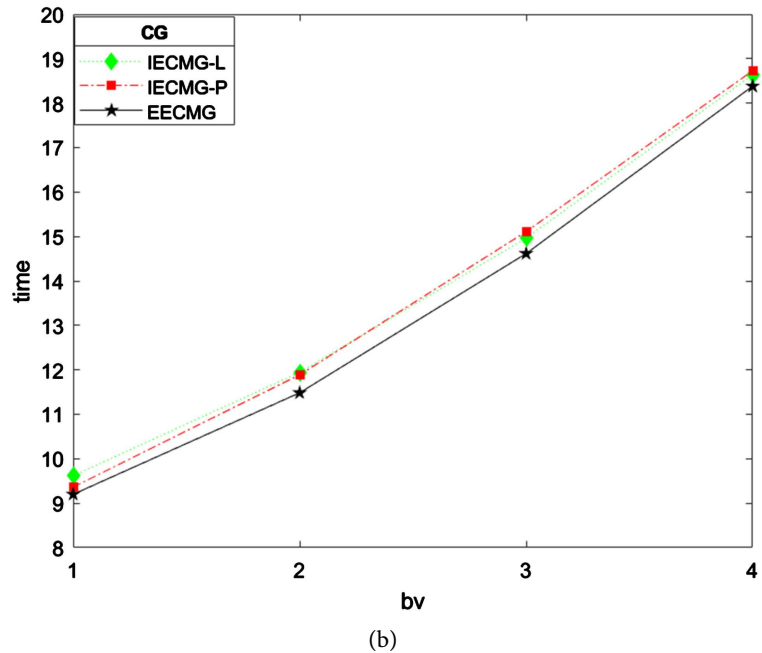

(b)

**Figure 6.** PSNR and Time of IECMG (CG) and EECMG (CG). (a) PSNR; (b) Time.

**Table 3.** Numerical results of IECMG (CG) and EECMG (CG).

| Blur and noise levels | IECMG-L(CG) | | IECMG-P(CG) | | EECMG (CG) | |
|---|---|---|---|---|---|---|
| | PSNR | Time | PSNR | Time | PSNR | Time |
| $b_1v_1$ | 46.8938 | 9.617 | 47.5923 | 9.362 | 54.2563 | 9.200 |
| $b_2v_2$ | 34.9539 | 11.943 | 35.4412 | 11.888 | 37.2658 | 11.484 |
| $b_3v_3$ | 28.9261 | 14.964 | 29.2944 | 15.109 | 31.0168 | 14.616 |
| $b_4v_4$ | 25.4300 | 18.636 | 25.7826 | 18.727 | 27.0067 | 18.381 |

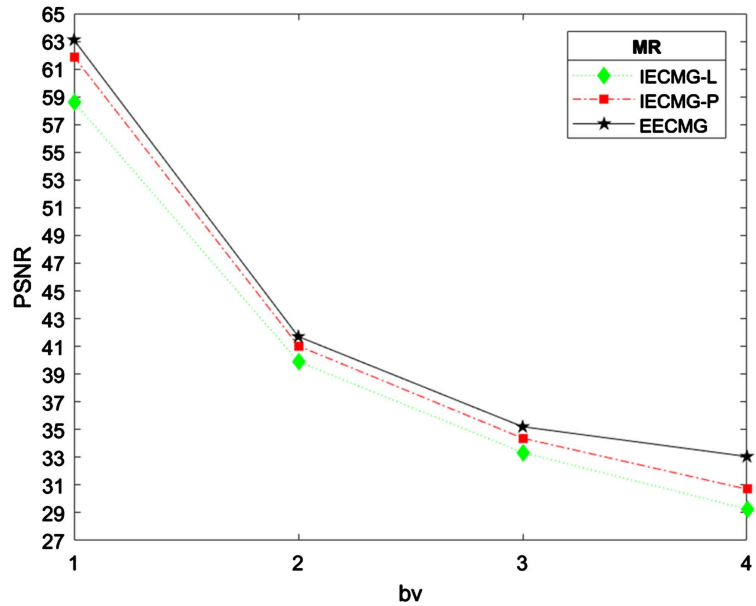

(a)

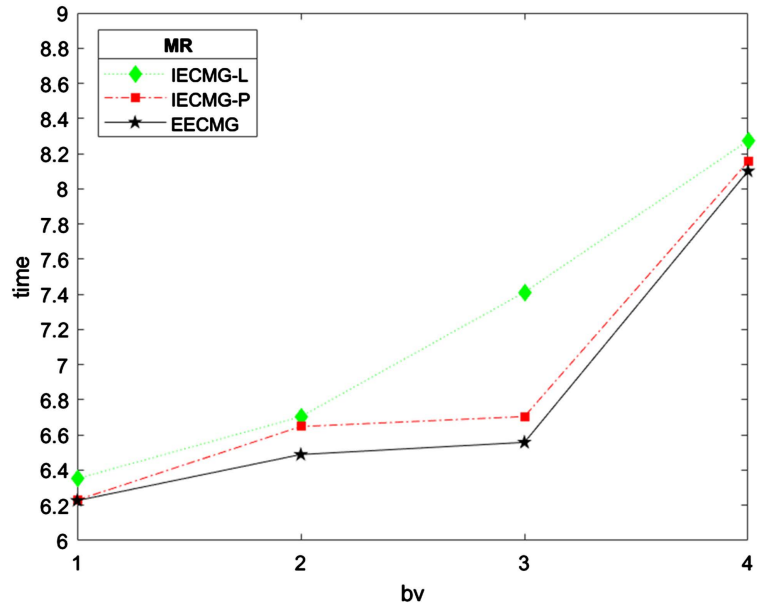

(b)

**Figure 7.** PSNR and TIME of IECMG (MR) and EECMG (MR). (a) PSNR; (b) Time.

**Table 4.** Numerical results of IECMG (MR) and EECMG (MR).

| Blur and noise levels | IECMG-L(MR) | | IECMG-P(MR) | | EECMG (MR) | |
|---|---|---|---|---|---|---|
| | PSNR | Time | PSNR | Time | PSNR | Time |
| $b_1 v_1$ | 58.6189 | 6.352 | 61.8852 | 6.230 | 63.0973 | 6.227 |
| $b_2 v_2$ | 39.9041 | 6.704 | 41.0350 | 6.648 | 41.7086 | 6.489 |
| $b_3 v_3$ | 33.3558 | 7.412 | 34.3789 | 6.705 | 35.2055 | 6.358 |
| $b_4 v_4$ | 29.2522 | 8.274 | 30.7030 | 8.155 | 33.0595 | 8.101 |

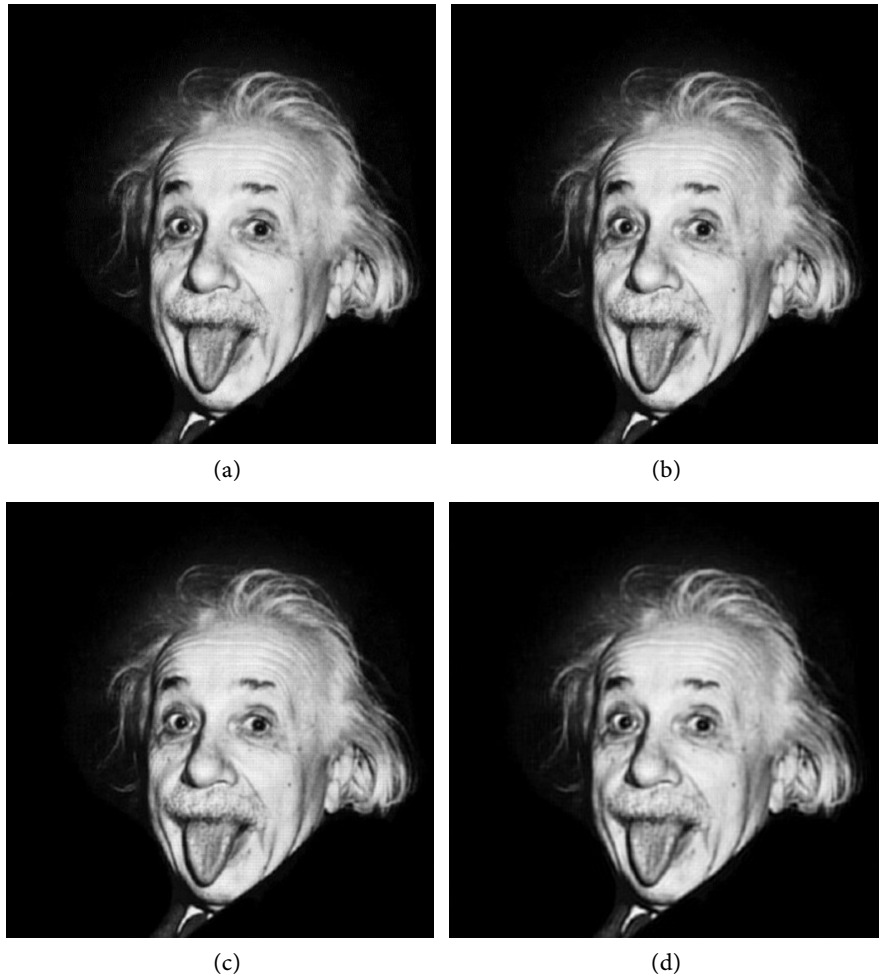

Figure 8. Recovery image of EECMG (MR) for four different blur and noise levels. (a) EECMG, $b_1v_1$, PSNR = 63.0973; (b) EECMG, $b_2v_2$, PSNR = 41.7086; (c) EECMG, $b_3v_3$, PSNR = 35.2055; (d) EECMG, $b_4v_4$, PSNR = 33.0595.

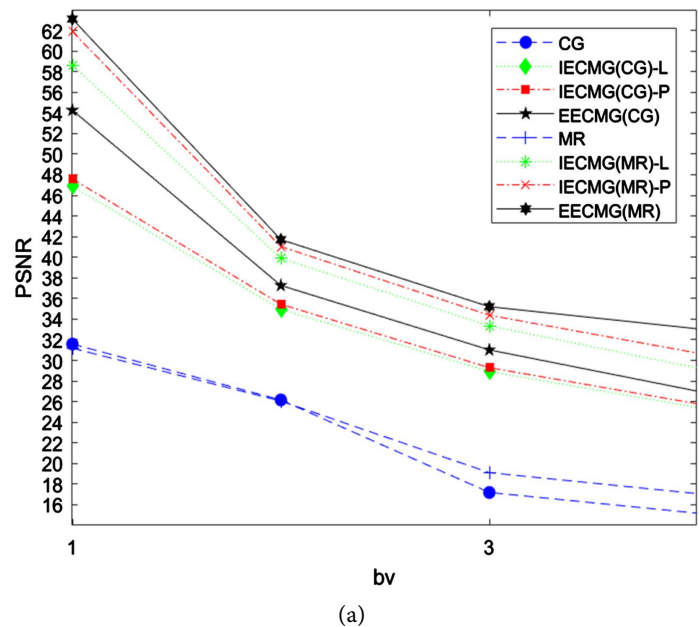

(a)

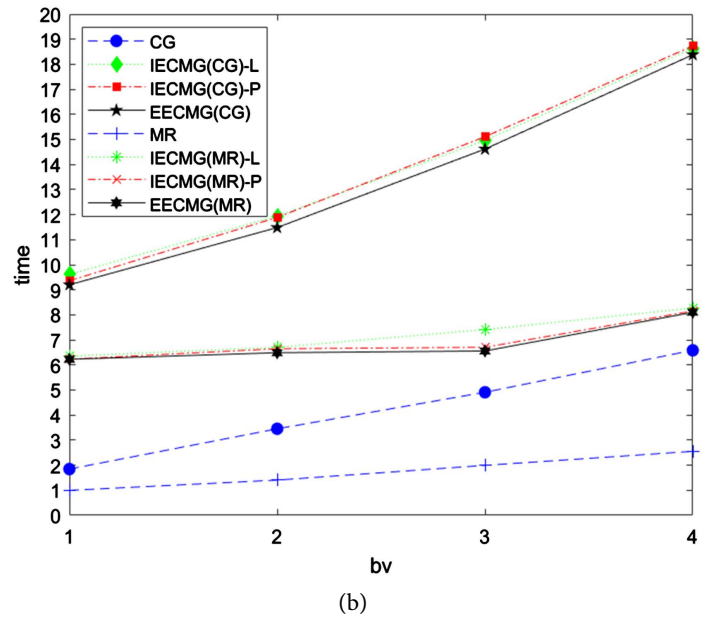

(b)

**Figure 9.** Numerical results of CG, MR, IECMG and EECMG. (a) PSNR; (b) Time.

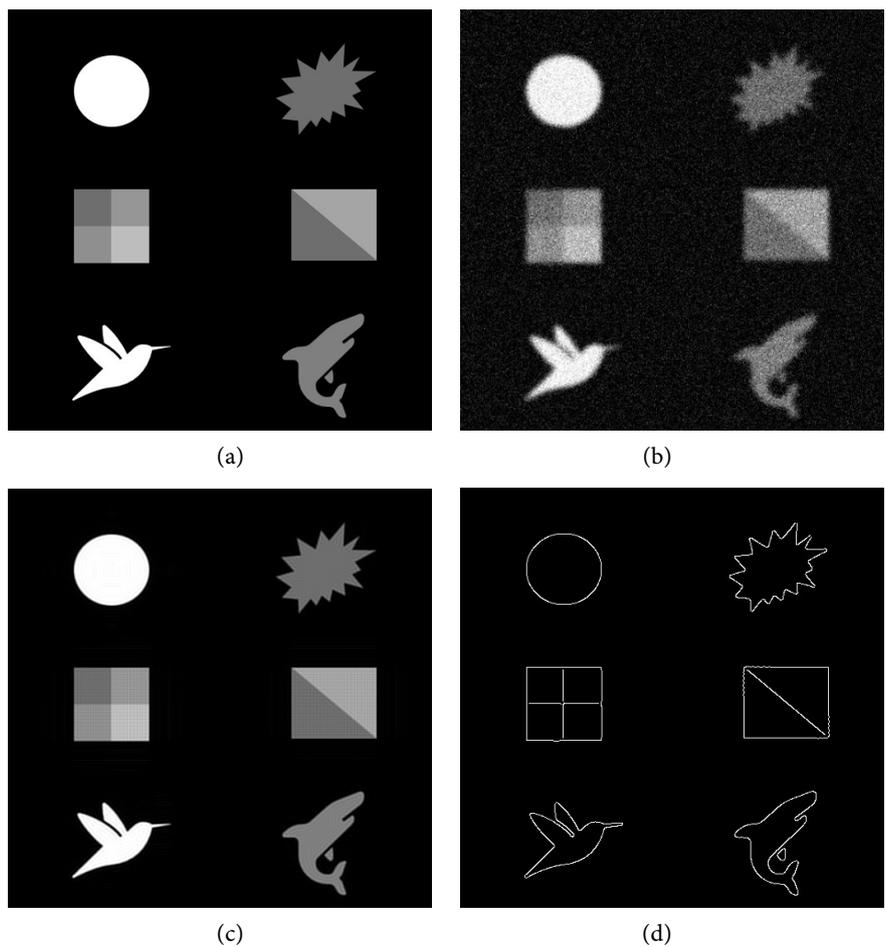

**Figure 10.** Edge preserving effect of EECMG (MR). (a) Original image; (b) Blur and noise image with $b_3v_3$; (c) Recovery image; (d) Recovery of image edge.

## 6. Conclusions

Combining the advantages of the new extrapolation method and the economical cascadic multigrid method, this paper proposes a new extrapolation economical cascadic multigrid method for image restoration, which not only improves the computational efficiency but also ensures the good quality of restoration. Numerical experiments show that,

1) The edge preserving denoising operator can achieve the effect of preserving image edges and removing noise;

2) The weighted least squares approximation smooths the image locally and reduces the staircase effect;

3) The new extrapolation formula and quadratic interpolation construct a nonlinear continuation operator, which can provide more accurate pixel values for the fine grid level and ensure good recovery quality;

4) The iteration number selection method controls the iteration number of each grid level, which greatly reduces the computational work.

Although this paper puts forward several kinds of algorithms for linear model of image restoration problem and has made some research achievements, there are still many works worth further research, such as:

1) Whether the algorithm proposed in this paper can also have good restoration effect on other image restoration models?

2) Whether the algorithm in this paper can restore large-size images?

3) Can the algorithm in this paper be applied to the restoration of three-dimensional images?

## Acknowledgements

This work was supported by Natural Science Foundation of China (11661027), the Guangxi Natural Science Foundation (2020GXNSFAA159143), and National Natural Science Foundation of China (12161027).

## Conflicts of Interest

The authors declare no conflicts of interest regarding the publication of this paper.